\documentclass[11pt]{article}
\usepackage{cite}
\usepackage{mathrsfs}
\usepackage{amsfonts}
\usepackage{amsmath}
\usepackage{amsfonts,amssymb,color}
\usepackage{dsfont}
\usepackage{curves}
\usepackage{mathrsfs}
\usepackage{pifont}
\usepackage{amssymb}
\usepackage{latexsym,amsmath,amssymb,amsfonts,epsfig,graphicx,cite,psfrag}
\usepackage{eepic,color,colordvi,amscd}
\usepackage{enumitem}

\newtheorem{theorem}{Theorem}[section]

\newtheorem{lemma}{Lemma}[section]

\newtheorem{claim}{Claim}[section]
\newtheorem{conjecture}{Conjecture}[section]

\newcommand{\qed}{\hfill\rule{0.5em}{0.809em}}

\def\emptyset{\mbox{{\rm \O}}}
\def\bar{\overline}

\renewcommand{\baselinestretch}{1.2}
\def\qed{\hfill \rule{4pt}{7pt}}

\def\pf{\noindent {\it Proof. }}

\begin{document}
	
	\title{Coloring of some $(P_2\cup P_4)$-free graphs}
		\author{Ran Chen$^{1,}$\footnote{Email: 1918549795@qq.com},  \; Xiaowen Zhang$^{2,}$\footnote{Corresponding author. Email: xiaowzhang0128@126.com}\\\\
		\small $^1$Institute of Mathematics, School of Mathematical Sciences\\
		\small Nanjing Normal University, 1 Wenyuan Road,  Nanjing, 210023,  China\\
		\small $^2$Department  of Mathematics\\
		\small East China Normal University, Shanghai, 200241, China
	}
	\date{}
	
	\maketitle
	
	\begin{abstract}
		We denote a path on $t$ vertices as $P_t$ and a cycle on $t$ vertices as $C_t$. For two vertex-disjoint graphs $G_1$ and $G_2$, the {\em union} $G_1\cup G_2$ is the graph with $V(G_1\cup G_2)=V(G_1)\cup V(G_2)$ and $E(G_1\cup G_2)=E(G_1)\cup E(G_2)$.   A {\em diamond} (resp. {\em gem}) is a graph consisting of a $P_3$ (resp. $P_4$) and a new vertex adjacent to all vertices of the $P_3$ (resp. $P_4$), and a {\em butterfly} is a graph consisting of two triangles that share one vertex. In this paper, we show that $\chi(G)\le 3\omega(G)-2$ if $G$ is a ($P_2\cup P_4$, gem)-free graph, $\chi(G)\le \frac{\omega(G)^2+3\omega(G)-2}{2}$ if $G$ is a ($P_2\cup P_4$, butterfly)-free graph. We also study the class of ($P_2\cup P_4$, diamond)-free graphs, and show that, for such a graph $G$,  $\chi(G)\leq4$ if $\omega(G)=2$, $\chi(G)\leq7$ if $\omega(G)=3$, $\chi(G)\leq9$ if $\omega(G)=4$, and $\chi(G)\leq2\omega(G)-1$ if $\omega(G)\ge 5$. Moreover, we prove that $G$ is perfect if $G$ is ($P_2\cup P_4$, diamond, $C_5$)-free with $\omega(G)\geq5$.
		
\begin{flushleft}
	{\em Key words and phrases:} ($P_2\cup P_4$)-free graphs, chromatic number, clique number\\
	{\em AMS 2000 Subject Classifications:}  05C15, 05C75\\
\end{flushleft}

	\end{abstract}

\renewcommand{\baselinestretch}{1.8}

\section{Introduction}\label{111}

All graphs considered in this paper are finite and simple. We follow \cite{BM08} for undefined notations and terminology.  We use $P_k$ and $C_k$ to denote a {\em path} and a {\em cycle} on $k$ vertices respectively. Let $G$ be a graph with vertex set $V(G)$ and edge set $E(G)$, and let $X\subseteq V(G)$. We use $G[X]$ to denote the subgraph of $G$ induced by $X$, and call $X$ a {\em clique} (resp. {\em stable set}) if $G[X]$ is a complete graph (resp. has no edges). The {\em clique number} of $G$, denoted by $\omega(G)$, is the maximum size taken over all cliques of $G$. 

We say that a graph $G$ contains a graph $H$ if $H$ is isomorphic to an induced subgraph of $G$, and say that $G$ is $H$-{\em free} if it does not contain $H$.
For a family $\{H_1,H_2,\cdots\}$ of graphs, $G$ is $(H_1, H_2,\cdots)$-free if $G$ is $H$-free for every $H\in \{H_1,H_2,\cdots\}$. For two vertex-disjoint graphs $G_1$ and $G_2$, the {\em union} $G_1\cup G_2$ is the graph with $V(G_1\cup G_2)=V(G_1)\cup V(G_2)$ and $E(G_1\cup G_2)=E(G_1)\cup E(G_2)$, and $G_1+G_2$ is the graph with $V(G_1+G_2)=V(G_1)\cup V(G_2)$ and $E(G_1+G_2)=E(G_1)\cup E(G_2)\cup \{uv~|~u\in V(G_1), v\in V(G_2)\}$.

A {\em k-coloring} of a graph $G=(V,E)$ is a mapping $f$: $V \rightarrow \{1,2,...,k\}$ such that $f(u)\neq f(v)$ whenever $uv\in E$. We say that $G$ is {\em k-colorable} if $G$ admits a $k$-coloring. The {\em chromatic number} of $G$, denoted by $\chi(G)$, is the smallest positive integer $k$ such that $G$ is $k$-colorable. 

The concept of binding functions was introduced by Gy\'{a}rf\'{a}s \cite{G75} in 1975. A family $\cal{G}$ of graphs is said to be {\em $\chi$-bounded} if there exists a function $f$ such that for every graph $G\in \cal{G}$, $\chi(G)\le f(\omega(G))$. If such a function does exist to $\cal{G}$, then the function $f$ is called a {\em binding function} for $\cal{G}$. 
 
An induced cycle of length at least 4 is called a {\em hole}, and
its complement is called an {\em antihole}. A hole or antihole is odd or
even if it has odd or even number of vertices. A graph $G$ is said to be {\em perfect} if $\chi(H)=\omega(H)$ for every induced subgraph $H$ of $G$. The famous Strong Perfect Graph Theorem \cite{CRSR06} states that a graph is perfect if and only if it is (odd hole, odd antihole)-free. Then Gy\'{a}rf\'{a}s \cite{G75}, and Sumner \cite{S81} independently, posed the following conjecture.
 
\begin{conjecture}\label{forest}{\em \cite{G75,S81}.}
	 For every forest $T$, the class of T-free graphs is $\chi$-bounded.
\end{conjecture}

 Gy\'{a}rf\'{a}s \cite{G87} proved the conjecture for $T=P_t$: every $P_t$-free graph $G$ has $\chi(G)\le (t-1)^{\omega(G)-1}$. Note that this binding function is exponential in $\omega(G)$. It's natural to ask the following question: Is it possible to improve the exponential bound for $P_t$-free graphs to a polynomial bound \cite{S16}. 

\begin{figure}[htbp]
	\begin{center}
	\includegraphics[width=17cm]{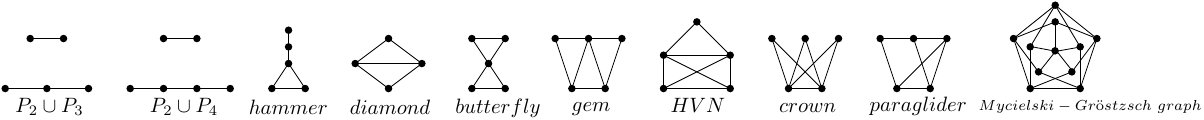}
	\end{center}
	\vskip -25pt
	\caption{Illustration of some forbidden configurations.}
	\label{fig-1}
\end{figure}

In \cite{BC18}, Bharathi {\em et al} showed that the binding function for $(P_2\cup P_3)$-free graphs is $f(\omega)=\frac{1}{6}\omega(\omega+1)(\omega+2)$, which represents the best bound so far. In \cite{PFR22} Prashant {\em et al}. obtained linear binding functions for the class of ($P_2\cup P_3$, $(K_1\cup K_2)+K_p$)-free graphs and ($P_2\cup P_3$, $2K_1+K_p$)-free graphs. In addition, they gave tight binding functions for ($P_2\cup P_3$, HVN)-free graphs and for ($P_2\cup P_3$, diamond)-free graphs. Wu and Xu \cite{WX23} proved that $\chi(G)\le\frac{1}{2}\omega^2(G)+\frac{3}{2}\omega(G)+1$ if $G$ is ($P_2\cup P_3$, crown)-free, Char and Karthick \cite{CK22} proved that $\chi(G)\le$ max $\{\omega(G)+3, \lfloor \frac{3\omega(G)}{2} \rfloor-1\}$ if $G$ is a ($P_2\cup P_3$, paraglider)-free graph with $\omega(G)\ge3$, Prashant {\em et al.} \cite{PAM23} proved that $\chi(G)\le2\omega(G)$ if $G$ is ($P_2\cup P_3$, gem)-free, and the bound was recently improved by Char and Karthick \cite{CK24} to $\lceil\frac{(5\omega(G)-1)}{4}\rceil$ for $\omega(G)\geq4$. In \cite{LLw2023}, \cite{LLW2023'} and \cite{LWL2023}, Li {\em et al} gave a linear binding function for $(P_2\cup P_3, H)$-free graphs, where $H\in$\{house, $K_1+C_4$, kite\}, and they also gave a binding function for $(P_2\cup P_3$, hammer)-free and $(P_2\cup P_3,C_5)$-free graphs. Very recently, Xu and Zhang \cite{XZ24} also proved $\chi(G)\le 4$ if $G$ is ($P_2\cup P_3, K_4$)-free graph with the distance of any two triangles of $G$ is at least $1$. See Figure~\ref{fig-1} for the illustration of some forbidden configurations.
 
In this paper, we focus on $(P_2\cup P_4)$-free graphs, which is a superclass of ($P_2\cup P_3$)-free graphs. Thus far, the best known binding function for $(P_2\cup P_4)$-free graphs is the same as that for $P_2\cup P_3$-free graphs: $f(\omega)=\frac{1}{6}\omega(\omega+1)(\omega+2)$ in \cite{BC18}. In \cite{TKM18}, Karthick and Mishra showed that $\chi(G)\leq6$ if $G$ is ($P_2\cup P_3$, diamond, $K_4$)-free, and this bound is tight because the graph $H$ which is the complement of the Schl\"afli graph (see www.win.tue.nl/~aeb/graphs/Schlaefli.html)satisfies $\chi(H)=6$. They also showed that $\chi(G)\leq6$ if $G$ is ($P_6$, diamond, $K_4$)-free. In \cite{BC18}, Bharathi and Choudum showed that $G$ is perfect if $G$ is ($P_2\cup P_3$, diamond)-free with $\omega(G)\geq5$. 

Let $G$ be a ($P_2\cup P_4$, diamond)-free graph. We prove that $\chi(G)\leq4$ if $\omega(G)=2$ (this bound is tight because of the Mycielski-Gr\"ostzsch graph as shown in Figure~\ref{fig-1}), $\chi(G)\le 7$ if $\omega(G)=3$, $\chi(G)\le 9$ if $\omega(G)=4$ and $\chi(G)\le 2\omega(G)-1$ if $\omega(G)\geq5$. We also prove that $G$ is perfect if $G$ is a ($P_2\cup P_4$, diamond, $C_5$)-free graph with $\omega(G)\ge 5$. Moreover, we prove that $\chi(G)\le 3\omega(G)-2$ if $G$ is a ($P_2\cup P_4$, gem)-free graph, and $\chi(G)\le \frac{\omega(G)^2+3\omega(G)-2}{2}$ if $G$ is a ($P_2\cup P_4$, butterfly)-free graph. Notice that the class of $2K_2$-free graphs does not admit a linear binding function \cite{BRSV19}, and thus one can not expect a linear binding function for the class of $(P_2\cup P_4,$ butterfly)-free graphs.

\begin{theorem}\label{gem}
	Let $G$ be a ($P_2\cup P_4$, gem)-free graph. Then $\chi(G)\le 3\omega(G)-2$.
\end{theorem}

\begin{theorem}\label{butterfly}
	Let $G$ be a ($P_2\cup P_4$, butterfly)-free graph. Then $\chi(G)\le \frac{\omega(G)^2+3\omega(G)-2}{2}$.
\end{theorem}

\begin{theorem}\label{diamond}
	Let $G$ be a ($P_2\cup P_4$, diamond)-free graph. Then:
	\[
\chi(G) \leq 
\left\{
\begin{array}{ll}
4 & \text{if } \omega(G) = 2 \\
7 & \text{if } \omega(G) = 3 \\
9 & \text{if } \omega(G) = 4\\
2\omega(G)-1 & \text{if } \omega(G)\ge 5.
\end{array}
\right.
\]

\end{theorem}
	
\begin{theorem}\label{C5}
	Let $G$ be a ($P_2\cup P_4$, diamond, $C_5$)-free graph with $\omega(G)\ge 5$. Then $G$ is perfect.
\end{theorem}

We will prove Theorem~\ref{gem} in Section~\ref{333}, prove Theorem~\ref{butterfly} in Section~\ref{444}, prove Theorems~\ref{diamond} in Section~\ref{555}, and prove Theorem~\ref{C5} in Section~\ref{666}.

\section{Preliminaries}	\label{222}
	
Let $G$ be a graph. For $v\in V(G)$, let $N_G(v)$ be the set of vertices adjacent to $v$, $d_G(v)=|N_G(v)|$, and $M_G(v)=V(G)\setminus (N_G(v)\cup\{v\})$. Let $N_G(X)=\{u\in V(G)\setminus X\;|\; u$ has a neighbor in $X\}$, and $M_G(X)=V(G)\setminus (X\cup N_G(X))$. If it does not cause any confusion, we usually omit the subscript $G$ and simply write $N(v)$, $d(v)$, $M(v)$, $N(X)$ and $M(X)$.

 For $S, T\subseteq V(G)$, let $N_T(S)=N(S)\cap T$ and $M_T(S)=T\setminus (N_T(S)\cup S)$. We use $[S,T]$ to denote the set of edges between $S$ and $T$. We say that $X$ is complete (resp. anticomplete) to $Y$  if $|[X,Y]|=|X||Y|$ (resp. $[X,Y]=\emptyset$), and we call $[X,Y]$ is complete if $X$ is complete to $Y$. For $u, v\in V(G)$, we simply write $u\sim v$ if $uv\in E(G)$, and write $u\not\sim v$ if $uv\not\in E(G)$.	
	
Throughout this paper we partition $V(G)$ into several vertex sets, which was initially defined by S. Wagon in \cite{W80} and later modified by A. P. Bharathi {\em et al} in \cite{BC18}.
Let $A$ be a maximum clique in $G$ with vertices $v_1, v_2 ,..., v_{\omega(G)}$. We define the sets $C_{i,j}$ for pairs of vertices $v_i, v_j$ of $A$ as follows: let $<_L$ denote the lexicographic ordering on the set $\{(i, j) \mid 1 \le i < j \le \omega(G)\}$, and for all $(i,j)$ in that set let $$C_{i,j}=\{v\in V(G)\setminus A \mid v\not\sim v_i~\text{and}~ v\not\sim v_j\}\setminus \bigcup_{(i',j')<_L (i,j)}C_{i',j'}.$$
There are $\binom{\omega(G)}{2}$ number of $C_{i,j}$ sets and these are pairwise disjoint. And every vertex in $C_{i,j}$ is adjacent to $v_k$ of $A$, where $1\le k<j, k\neq i$. So each vertex $v\in V(G)\setminus(A\cup (\bigcup_{i,j}C_{i,j}))$ is adjacent to $|A|-1$ vertices of $A$. Then we define the following sets. For $v_a\in A$, let
$$I_a=\{v\in V(G) \mid v\sim x, \forall x\in A\setminus\{v_a\}~\text{and}~v\not\sim v_a\}.$$
By definition, we have that $$V(G)=A\cup (\bigcup_{i,j}C_{i,j})\cup (\bigcup_{v_a\in A}I_a).$$

In this paper, we suppose $G$ is always a ($P_2\cup P_4$)-free graphs, hence we demonstrate some properties of $C_{i,j}$ and $I_a$. To prove these properties, we need the following Lemma.

\begin{lemma}{\em \cite{SD74}}\label{P_4}
	Every $P_4$-free graph is perfect. 
\end{lemma}

\begin{claim}\label{perfect}
	Each $G[C_{i,j}]$ is $P_4$-free and hence it is perfect, and each $I_a$ is a stable set. 
\end{claim}

\pf If there exists an induced $P_4=xyzw$ in $G[C_{i,j}]$, then $G[\{v_i,v_j\}\cup \{x,y,z,w\}]$ is an induced $P_2\cup P_4$, a contradiction. Suppose $G[I_a]$ contains an edge $xy$, but then $A\cup \{x,y\}-\{v_a\}$ is a clique of size $\omega(G)+1$. Both are contradictions. By Lemma~\ref{P_4}, Claim~\ref{perfect} holds.\qed

\section{Proof of Theorem~\ref{gem}}\label{333}

Let $G$ be a ($P_2\cup P_4$, gem)-free graph. In this section, we will prove Theorem~\ref{gem}. Before that, we give some structural properties of $G$. All notations are defined the same as in Section~\ref{222}.

Let $M=(\bigcup_{k\ge 2}I_k)\cup (\bigcup_{i\ge 2}C_{i,j})\cup \{v_2,..., v_{\omega(G)}\}$ and $N=I_1\cup (\bigcup_{j\ge 3}C_{1,j})$. Clearly, we have that
\begin{equation}\label{eqa-1}
	\mbox{$[v_1,M]$ is complete}.
	\end{equation}
and
\begin{equation}\label{eqa-2}
	\mbox{$[v_2,N]$ is complete}. 
\end{equation}

By (\ref{eqa-1}) and (\ref{eqa-2}), and noting that $G$ is gem-free, we have that:
\begin{equation}\label{eqa-3}
	\mbox{$G[M]$ and $G[N]$ are both $P_4$-free and hence they are perfect graphs.}
\end{equation}	
	
By (\ref{eqa-1}), (\ref{eqa-2}) and (\ref{eqa-3}), $\chi(G[M])=\omega(G[M]) \le \omega(G)-1$ and $\chi(G[N])=\omega(G[N])\le \omega(G)-1$. Hence 
\begin{align*}
	\chi(G)&\le \chi(G[M])+\chi(G[N])+\chi(G[C_{1,2}])\\
	& \le (\omega(G)-1)+(\omega(G)-1)+\omega(G)\\
	& = 3\omega(G)-2.
\end{align*}

It completes the proof of Theorem~\ref{gem}.\qed

\section{Proof of Theorem~\ref{butterfly}}\label{444}
	
In this section, we primarily study ($P_2\cup P_4$, butterfly)-free graphs. Now let $G$ be a ($P_2\cup P_4$, butterfly)-free graph, and all notations are defined the same as in Section~\ref{222}.  Before proving Theorem~\ref{butterfly}, we present an additional property of $C_{i,j}$.
\begin{equation}
	\mbox{$C_{i,j}$ is a stable set  where $j\ge 3$.}
\end{equation}

Suppose not, and let $xy\in E(G[C_{i,j}])$ where $j\ge 3$. By definition of $C_{i,j}$, each vertex in $C_{i,j}$ is adjacent to a vertex $z\in \{v_1, v_2\}\subseteq A$. Then $\{x,y,z,v_i,v_j\}$ induces a butterfly, a contradiction.

Since $\omega(C_{1,2})\le \omega(G)$ and $G[C_{1,2}]$ is a perfect graph by Claim~\ref{perfect}, we have that $\chi(G[C_{1,2}])=\omega(G[C_{1,2}])\le \omega(G)$. Then
\begin{align*}
	\chi(G[C])&\le \chi(G[C_{1,2}])+\chi(G[C_{1,3}])+\ldots+\chi(G[C_{\omega(G)-1, \omega(G)}])\\
	& \le \omega(G)+\binom{\omega(G)}{2}-1\\
	& = \frac{\omega(G)^2+\omega(G)-2}{2}
\end{align*}

Moreover, $I_a$ is a stable set by Claim~\ref{P_4}, and it is anticomplete to $v_a$ in $A$. Consequently, we can assign the same color to $G[I_a]$ as that of $v_a$.

Therefore 
\begin{align*}
	\chi(G)&\le \chi(G[A])+\chi(G[C])\\
	& \le \omega(G)+ \frac{\omega(G)^2+\omega(G)-2}{2}\\
	& = \frac{\omega(G)^2+3\omega(G)-2}{2}
\end{align*}	

It completes the proof of Theorem~\ref{butterfly}.\qed

\section{Proof of Theorem~\ref{diamond}} \label{555} 

In this section, we consider ($P_2\cup P_4$, diamond)-free graphs. Let $G$ be a ($P_2\cup P_4$, diamond)-free graph,  and all notations are defined the same as in Section~\ref{222}. If $\omega(G)=1$, then obviously the chromatic number of $G$ is $1$. So in the following, we may suppose $\omega(G)\ge 2$. 

\begin{claim}\label{emptyset}
	$C_{i,j}=\emptyset$ for $j\ge 4$.
\end{claim}

\pf Suppose, to the contrary, that $x\in C_{i,j}$ for some $j\ge 4$. Then according to the definition of $C_{i,j}$, there exist two distinct vertices $p,q\in\{v_1,v_2,v_3\}\subseteq A$ such that $x\sim p$ and $x\sim q$. But now $\{x,v_j,p,q\}$ induces a diamond, a contradiction.\qed

\begin{claim}\label{I}
	For $1\le a\le \omega(G)$, $I_a=\emptyset$ if $\omega(G)\ge 3$, and it is a stable set if $\omega(G)=2$.  
\end{claim}
\pf If $\omega(G)\ge 3$, and $x\in I_a$. By the definition of $I_a$, $[x,A\setminus \{v_a\}]$ is complete. Then there exists two distinct vertices $p,q\in A\setminus \{v_a\}$ such that $x\sim p$ and $x\sim q$. Hence $\{x,p,q,v_a\}$ induces a diamond, a contradiction. If $\omega(G)=2$, then the assertion follows by Claim~\ref{P_4}.
\qed

\medskip
By Claims~\ref{emptyset} and \ref{I}, we have that 
$$V(G)=A\cup C_{1,2}\cup I_1\cup I_2,~\text{if}~\omega(G)=2;$$
$$V(G)=A\cup C_{1,2}\cup C_{1,3}\cup C_{2,3},~\text{if}~\omega(G)\ge 3.$$

Moreover, by the definitions of $C_{1,3}$ and $C_{2,3}$,
\begin{equation}\label{complete}
	\mbox{ $v_2$ is complete to $C_{1,3}$ and $v_1$ is complete to $C_{2,3}$.}
\end{equation}


By Claim~\ref{perfect} and (\ref{complete}), we have that
\begin{equation}\label{chi}
	\mbox{$\chi(C_{12})=\omega(G[C_{1,2}])$, ~$\chi(G[C_{1,3}])=\omega(G[C_{1,3}])\le \omega(G)-1$~and~$\chi(G[C_{2,3}])=\omega(G[C_{2,3}])\le \omega(G)-1$.}
\end{equation}

\begin{claim}\label{C3}
	$[C_{1,3}, A\setminus \{v_2\}]=\emptyset$, and $[C_{2,3}, A\setminus \{v_1\}]=\emptyset$.
\end{claim}

\pf If there exists an edge $xv_i\in [C_{1,3}, A\setminus \{v_2\}]$, where $x\in C_{1,3}$, $v_i\in A\setminus\{v_2\}\}$, then $\{x,v_i,v_1,v_2\}$ induces a diamond, a contradiction. Similarly, $[C_{2,3}, A\setminus \{v_1\}]=\emptyset$.
\qed

\begin{claim}\label{C12}
	Each vertex of $C_{1,2}$ is adjacent to at most one vertex of $A$.
\end{claim}
\pf This claim is obvious if $\omega(G)\leq3$. Now, assume that $\omega(G)\ge 4$. Suppose some vertex $x\in C_{1,2}$ is adjacent to two distinct vertices $p$ and $q$ in $A
\setminus \{v_1,v_2\}$. Then $\{x,p,q,v_1\}$ induces a diamond, a contradiction.
\qed
\medskip

We now prove Theorem~\ref{diamond} depending on different values of $\omega(G)$.

\begin{itemize}
    \item $\omega(G)=2$, thus $A=\{v_1, v_2\}$.
   \end{itemize}
   
 We have that $V(G)=A\cup C_{1,2}\cup I_1\cup I_2$. According to Claim~\ref{I} and (\ref{chi}), we can color $A$ with 1,2, color $I_1$ with 1, color $I_2$ with 2, and color $C_{1,2}$ with 3,4, thus $\chi(G)\le 4$.
 
\begin{itemize}
    \item $\omega(G)=3$, thus $A=\{v_1, v_2,v_3\}$.
   \end{itemize} 
 
 We have that $V(G)=A\cup C_{1,2}\cup C_{1,3}\cup C_{2,3}$ when $\omega(G)=3$. Moreover $\chi(G[C_{1,3}])\le 2$, $\chi(G[C_{2,3}])\le 2$ and $\chi(G[C_{1,2}])\le 3$ by (\ref{chi}). Now, we can color $A$ with 1,2,3, color $C_{1,2}$ with 1,2,4, color $C_{2,3}$ with 3,5 and color $C_{1,3}$ with 6,7. Therefore we have that $\chi(G)\le 7$.

\begin{itemize}
    \item $\omega(G)=4$, thus $A=\{v_1, v_2, v_3, v_4\}$.
   \end{itemize} 

Recall that $V(G)=A\cup C_{1,2}\cup C_{1,3}\cup C_{2,3}$. Let $Q$ be a maximum clique of $G[C_{1,2}]$. We have that $\chi(G[C_{1,2}])\leq|Q|\leq4$, and $\max\{\chi(G[C_{1,3}]), \chi(G[C_{2,3}])\}\leq3$ by (\ref{chi}). 

First, we consider the case where $|Q|\le 3$, then $\chi(G[C_{1,2}])\leq3$. According to Claim~\ref{C3}, we can color $A$ with 1,2,3,4, color $C_{1,2}$ with 1,2,5, color $C_{1,3}$ with 3,4,6 and color $C_{2,3}$ with 7,8,9. And so, $\chi(G)\leq9$ if $|Q|\leq3$.

Second, we consider the case where $|Q|=4$. Suppose that $C_{1,3}=\emptyset$. Then we can color $A$ with 1,2,3,4, color $C_{1,2}$ with 1,2,5,6 and color $C_{2,3}$ with 3,4,7. Hence $\chi(G)\le 7$. Similarly, we have that $\chi(G)\le 7$ if $C_{2,3}=\emptyset$. Now, suppose that $C_{1,3}\neq \emptyset$ and $C_{2,3}\neq \emptyset$. We prove that 
\begin{equation}\label{4-1}
	\mbox{$C_{1,3}$ is complete to $C_{2,3}$.}
\end{equation}

Let $x\in C_{1,3}$ and $y\in C_{2,3}$ such that $x\not\sim y$, and let $Q=\{u_1,u_2,u_3,u_4\}$. By definition, $x\sim v_2$ and $y\sim v_1$. Since $|Q|=\omega(G)=4$ and $G$ is diamond-free, $|N_Q(x)|\le 1$ and $|N_Q(y)|\le 1$. Without loss of generality, we assume that $[\{u_1,u_2\}, \{x,y\}]=\emptyset$. Then $\{u_1,u_2\}\cup \{x,v_2,v_1,y\}$ induces a $P_2\cup P_4$, a contradiction. This proves (\ref{4-1}). 

Next we prove that 
\begin{equation}\label{4-2}
	\mbox{$C_{1,3}$ and $C_{2,3}$ are both stable sets.}
\end{equation}

Suppose that $C_{1,3}$ is not a stable set. Let $x_1x_2\in E(G[C_{1,3}])$ and $x_3\in C_{2,3}$. By (\ref{4-1}), $x_1\sim x_3$ and $x_2\sim x_3$. By (\ref{complete}) and Claim~\ref{C3}, $\{v_2,x_1,x_2,x_3\}$ induces a diamond, a contradiction. Similarly, $C_{2,3}$ is a stable set. 

Thus, we can color $A$ with 1,2,3,4, color $C_{1,2}$ with 1,2,5,6, color $C_{1,3}$ with 3 and color $C_{2,3}$ with 4 by (\ref{4-2}) and Claim~\ref{C3}. So, $\chi(G)\le 6$. Consequently, we have that $\chi(G)\leq9$ if $\omega(G)=4$.

\begin{itemize}
    \item $\omega(G)\ge 5$, thus $\{v_1, v_2, v_3, v_4, v_5\}\subseteq A$.
   \end{itemize} 

Now $V(G)=A\cup C_{1,2}\cup C_{1,3}\cup C_{2,3}$. Let $Q$ be a maximum clique of $G[C_{1,2}]$. By (\ref{chi}), $\chi(G[C_{1,2}])=|Q|\leq\omega(G)$, $\chi(G[C_{1,3}])=\omega(G[C_{1,3}])\leq\omega(G)-1$ and $\chi(G[C_{2,3}])=\omega(G[C_{2,3}])\leq\omega(G)-1$. Then we will consider the following four cases depending on the values of $|Q|$.

\medskip

{\noindent \bf  Case 1 $|Q|=1.$}

In this case, $C_{1,2}$ is a stable set. According to Claim~\ref{C3},  we can color $A$ with $1,2,\ldots,\omega(G)$, color $C_{1,2}$ with 2, color $C_{1,3}$ with $1,3,\ldots,\omega(G)$, and color $C_{2,3}$ with ($\omega(G)-1$) additional new colors. Hence $\chi(G)\le \omega(G)+\omega(G)-1=2\omega(G)-1$.

\medskip

{\noindent\bf Case 2 $|Q|=2.$}

Let $Q=\{v_1',v_2'\}$. For $i\in \{1,2\}$, we define $N_{i1}=\{x\in C_{i,3}\mid N_{Q}(x)=\{v_1'\}\}$, $N_{i2}=\{x\in C_{i,3}\mid N_{Q}(x)=\{v_2'\}\}$, $N_{i3}=\{x\in C_{i,3}\mid x~\text{is complete to}~Q\}$, and $N_{i4}=\{x\in C_{i,3}\mid [x, Q]=\emptyset\}$. By definition, we have that 
\begin{equation}\label{eq1}
	\mbox{$C_{1,3}=N_{11}\cup N_{12}\cup N_{13}\cup N_{14}$ and $C_{2,3}=N_{21}\cup N_{22}\cup N_{23}\cup N_{24}$}.
\end{equation}

For $i\in \{1,2\}$, the following properties holds.

{\noindent\bf (M1)}
	$N_{i1}$, $N_{i2}$ and $N_{i3}$ are stable sets.
	
	We may by symmetry set $i=1$. Suppose to the contrary that $x,x'\in N_{11}$ and $x\sim x'$. By the definition of $N_{11}$, $x\sim v_1'$ and $x'\sim v_1'$. Recall that $v_2$ is complete to $C_{1,3}$, we have that $v_2\sim x$ and $v_2\sim x'$. Then $\{v_1',x,x',v_2\}$ induces a diamond, a contradiction. Similarly, $N_{12}$ and $N_{1,3}$ are also stable sets. This proves (M1).
\medskip

{\noindent\bf (M2)}	$[N_{i1}\cup N_{i2}, N_{i3}]=\emptyset$.

By symmetry, set $i=1$. Suppose $y\in N_{11}$ and $y'\in N_{13}$ such that $y\sim y'$. By the definitions of $N_{11}$ and $N_{13}$, $y\sim v_1'$ and $y'\sim v_1'$. Since $v_2$ is complete to $\{y,y'\}$, there exists an induced diamond on $\{v_1',y,y',v_2\}$, a contradiction. Similarly, $[N_{12},N_{1,3}]=\emptyset$. So (M2) holds.
\medskip

{\noindent\bf (M3)} Either $N_{14}=\emptyset$ or $N_{24}=\emptyset$.

 Suppose, to the contrary, that $z\in N_{14}$ and $z'\in N_{24}$. We have that $z\sim v_2$ and $z'\sim v_1$. Since $[\{z,z'\}, Q]=\emptyset$, to avoid a $P_2\cup P_4$ on $\{v_1',v_2'\}\cup \{z,v_2,v_1,z'\}$, $z\sim z'$. Since $\omega(G)\ge 5$, there exists a vertex in $A\setminus \{v_1,v_2\}$, say $v_k$, that is not adjacent to either $v_1'$ or $v_2'$ by Claim~\ref{C12}. Now $\{v_1',v_2'\}\cup \{z,z',v_1,v_k\}$ induces a $P_2\cup P_4$, a contradiction. So, (M3) holds.
\medskip

According to (M3), without loss of generality, we may assume that $N_{24}=\emptyset$. By (M1) and (M2), we can color $A$ with $\{1,2,\ldots,\omega(G)\}$, color $N_{11}\cup N_{13}$ with 1, color $N_{21}\cup N_{23}$ with 2, color $N_{12}$ with 3, color $N_{22}$ with 4.

Now, we consider the coloring of $N_{14}$. By (\ref{complete}), we have that $C_{1,3}$ is complete to $v_2$, and thus $C_{1,3}$ is $P_3$-free because of the diamond-freeness of $G$. Notice that $C_{1,3}=N_{11}\cup N_{12}\cup N_{13}\cup N_{14}$ by (\ref{eq1}). Let $Q'$ be a component of $N_{14}$. Clearly, $Q'$ is a clique of size at most $\omega(G)-1$. If $|Q'|\leq\omega(G)-2$, then we may color $Q'$ with 5, 6,$\cdots$, $\omega(G)$ and 2 more additional colors by Claim~\ref{C3}. If $|Q'|=\omega(G)-1$, then $[N_{11}\cup N_{12}\cup N_{13},Q']=\emptyset$ as otherwise let $x\in N_{12}\cup N_{12}\cup N_{13}$ such that $x$ has a neighbor in $Q'$, then $x$ is complete to $Q'$ since $C_{1,3}$ is $P_3$-free; but now, $\{x,v_2\}\cup Q'$ is a clique of size $\omega(G)+1$, a contradiction. Now, we may color $Q'$ with 1, 3, 5, 6,$\cdots$, $\omega(G)$ and one more additional color. Therefore, we may color $N_{14}$ with 5, 6,$\cdots$, $\omega(G)$ and at most 2 more additional colors. 

Last, we may use 2 more additional new colors for $C_{1,2}$ since $\chi(G[C_{1,2}])=|Q|=2$. Therefore we have that $\chi(G)\le \omega(G)+2+2=\omega(G)+4\leq2\omega(G)-1$ because $\omega(G)\geq5$.

\medskip

{\noindent\bf Case 3 $|Q|=\omega(G).$}

In this case, we have the following property.
\begin{equation}\label{R1}
	\mbox{Either $C_{1,3}=\emptyset$ or $C_{2,3}=\emptyset$.}
\end{equation}

Suppose, to the contrary, that $x\in C_{1,3}$ and $y\in C_{2,3}$. Then $x\sim v_2$ and $y\sim v_1$ by (\ref{complete}). If $[x,Q]$ is complete, then we have an $(\omega(G)+1)$-clique. So, let $t_1$ be a nonneighbor of $x$ in $Q$. If $|N_Q(x)|\geq2$, then let $t_3,t_4\in N_Q(x)$, and $G[\{t_1,t_2,t_3,x\}]$ is a diamond, a contradiction.  It implies that $|N_Q(x)|\leq1$. Similarly, we also have that $|N_{Q}(y)|\le 1$. Moreover, $|Q|=\omega(G)\ge 5$, there exists two distinct vertices $u$ and $v$ in $Q$ such that $[\{u,v\}, \{x,y\}]=\emptyset$. To avoid a $P_2\cup P_4$ on $\{u,v\}\cup \{x,v_2,v_1,y\}$, we have that $x\sim y$. By Claim~\ref{C12}, there must be a vertex in $A\setminus \{v_1,v_2\}$, say $v_l$, such that $v_l\not\sim u$ and $v_l\not\sim v$ as $|A|\ge 5$. But now $\{u,v\}\cup \{x,y,v_1,v_l\}$ induces a $P_2\cup P_4$, which is a contradiction. This proves (\ref{R1}).

Without loss of generality, we may assume that $C_{2,3}=\emptyset$. In this case $$V(G)=A\cup C_{1,2}\cup C_{1,3}$$
Then by Claim~\ref{C3}, we can color $A$ with $\{1,2,\ldots, \omega(G)\}$, color $C_{1,3}$ with $\{1,3,\ldots,\omega(G)\}$, and color $C_{1,2}$ with 2 and ($\omega(G)-1$) additional new colors. Therefore we have that $\chi(G)\le \omega(G)+(\omega(G)-1)= 2\omega(G)-1$.

\medskip

{\noindent\bf Case 4 $3\le |Q|\le \omega(G)-1.$ }

For $i\in \{1,2\}$, we define $N_{i1}=\{x\in C_{i,3}\mid [x, Q]~\text{is complete}\}$, $N_{i2}=\{x\in C_{i,3}\mid |N_Q(x)|=1\}$ and $N_{i3}=\{x\in C_{i,3}\mid [N_{13},Q]=\emptyset\}$.
Suppose $C_{1,3}\neq N_{11}\cup N_{12}\cup N_{13}$. Let $x\in C_{1,3}\setminus (N_{11}\cup N_{12}\cup N_{13})$. Then there exist a vertex $y\in Q$ such that $x\not\sim y$ and two distinct vertices $u,v\in Q$ such that $x\sim u$ and $x\sim v$. Thus $\{x,u,v,y\}$ induces a diamond, which is a contradiction. Therefore
\begin{equation}
	C_{1,3}=N_{11}\cup N_{12}\cup N_{13}~\text{and}~C_{2,3}=N_{21}\cup N_{22}\cup N_{23}.
\end{equation}

The following properties hold.

{\noindent\bf (T1)}	$N_{11}$, $N_{12}$, $N_{21}$ and $N_{22}$ are all stable sets.

 Let $p\in Q$. Suppose $N_{11}$ contains an edge $p_1p_2$.  
 By the definition of $N_{11}$, $p_1\sim p$ and $p_2\sim p$. Recall that $v_2$ is complete to $C_{1,3}$ and $v_2$ is anticomplete to $C_{1,2}$, then $[\{p,p_1,p_2,v_2\}]$ is an induced diamond, a contradiction. With the same argument, $N_{21}$ is also a stable set. 

Suppose to the contrary that $N_{12}$ contains an edge $ab$. If $N_Q(a)=N_Q(b)=\{c\}$, then $\{c,a,b,v_2\}$ induces a diamond, which implies that $N_Q(a)\neq N_Q(b)$. Assume that $N_Q(a)=\{z\}$ and $N_Q(b)=\{z'\}$. Since $|Q|\ge 3$, $\exists z''\in Q\setminus\{z,z'\}$ such that $z''\not\sim a$ and $z''\not\sim b$. Moreover, by Claim~\ref{C12}, there exists a vertex in $A\setminus \{v_1,v_2\}$, say $v_m$, such that $[v_m,\{z,z''\}]=\emptyset$. Now $\{v_1,v_m\}\cup \{a,b,z,z''\}$ induces a $P_2\cup P_4$, a contradiction. Similarly, $N_{22}$ is a stable set. This proves (T1).
\medskip

{\noindent\bf (T2)}	$[N_{11},N_{12}]=\emptyset$ and $[N_{21},N_{22}]=\emptyset$.

Let $a_1\in N_{11}$ and $a_2\in N_{12}$ be two adjacent vertices. Suppose $N_Q(a_2)=\{u\}$. Then $a_1\sim u$ by the definition of $N_{11}$. Moreover, $v_2\not\sim u$ and $v_2$ is complete to $\{a_1,a_2\}$, then $\{u,a_1,a_2,v_2\}$ induces a diamond, a contradiction. Similarly, we can prove $[N_{21},N_{22}]=\emptyset$. This proves (T2).

\medskip
According to (T1) and (T2), we have that 
\begin{equation}\label{eqa4-1}
	\mbox{$N_{11}\cup N_{12}$ and $N_{21}\cup N_{22}$ are stable sets}.
\end{equation}

{\noindent\bf (T3)}	Either $N_{13}=\emptyset$ or $N_{23}=\emptyset$.

 Suppose that $u_1\in N_{13}$ and $u_2\in N_{23}$. We have that $u_1\sim v_2$ and $u_2\sim v_1$. Let $\{w_1,w_2\}\subseteq Q$. To avoid a $P_2\cup P_4$ on $\{w_1,w_2\}\cup \{u_1,v_2,v_1,u_2\}$, $u_1\sim u_2$ as $[\{u_1,u_2\}, Q]=\emptyset$. Moreover, $\omega(G)\ge 5$, there exists a vertex in $A\setminus \{v_1,v_2\}$, say $v_s$, such that $[v_s,\{w_1,w_2\}]=\emptyset$ by Claim~\ref{C12}. Now $\{w_1,w_2\}\cup \{u_1,u_2,v_1,v_s\}$ induces a $P_2\cup P_4$, which is a contradiction. This proves (T3).
\medskip

By (T3), without loss of generality, assume that $N_{23}=\emptyset$. Now 
$$V(G)=A\cup C_{1,2}\cup N_{11}\cup N_{12}\cup N_{13}\cup N_{21}\cup N_{22}$$
Then we also have the following property.

{\noindent\bf (T4)}	$[C_{1,2}, N_{13}]=\emptyset$.

 Suppose that $x_1\in N_{13}$ and $x_2\in C_{1,2}$ are two adjacent vertices. Since $Q$ is a maximum clique of $C_{1,2}$, there must exist a vertex $y_1$ such that $x_2\not\sim y_1$. Since $x_2$ has two neighbors, say $y_2$, $y_3$, in $Q$ leads to a diamond on $\{y_1,y_2,y_3,x\}$, we have that $|N_Q(x_2)|\leq1$. Moreover, since $|Q|\geq3$, there exists two distinct vertices $t_1$, $t_2$ in $Q$ such that $[x_2,\{t_1,t_2\}]=\emptyset$. Moreover $[\{v_1,v_2\},\{t_1,t_2,x_2\}]=\emptyset$, $x_1\sim v_2$ and $x_1\not\sim v_1$ by definition. Hence $\{t_1,t_2\}\cup \{x_2,x_1,v_2,v_1\}$ induces a $P_2\cup P_4$, a contradiction. This proves (T4).

We can now obtain a $(2\omega(G)-1)$-coloring of $G$. First, color $A$ with $\{1,2,\ldots,\omega(G)\}$. Next, color $N_{11}\cup N_{12}$ with 1, color $N_{21}\cup N_{22}$ with 2 by (\ref{eqa4-1}). Finally, color $C_{1,2}$ and $N_{13}$ with the same ($\omega(G)-1$) additional new colors by (T4) and (\ref{chi}). Then we have that $\chi(G)\le 2\omega(G)-1$.
\medskip

It completes the proof of Theorem~\ref{diamond}.\qed

\section{Proof of Theorem~\ref{C5}}\label{666}

In this section, we will prove Theorem~\ref{C5}. Let $G$ be a ($P_2\cup P_4,$ diamond, $C_5$)-free graph with $\omega(G)\geq5$. Clearly, each odd hole of length at least 9 contains a $P_2\cup P_4$, and each odd antihole of vertices number at least 7 contains a diamond. If $G$ is $C_7$-free, then $G$ is ($C_{2k+1}, \overline{C_{2k+1}}$)-free for all $k\ge 2$. Hence by the Strong Perfect Theorem~\cite{CRSR06}, $G$ is a perfect graph, and we are done. Therefore, we next show that $G$ is $C_7$-free. 

Suppose to its contrary that $G$ contains a $C_7$, say $C$. Notice that the Claims~\ref{emptyset}, \ref{I}, \ref{C3} and \ref{C12} proved in Section~\ref{555} hold. Also, in section~\ref{555}, we show that $V(G)=A\cup C_{1,2}\cup C_{1,3}\cup C_{2,3}$ whenever $\omega(G)\ge 5$. Recall that $A=\{v_1,\cdots,v_{\omega(G)}\}$. Since there are at most two vertices of $C$ can belong to $A$, there must exists an induced $P_4$ which is an induced subgraph of $C$, say $P$, with $V(P)=\{a,b,c,d\}$, such that $V(P)\subseteq C_{1,2}\cup C_{1,3}\cup C_{2,3}$. Since $C_{1,2}$ is $P_4$-free, we have that $|V(P)\cap C_{1,2}|\leq3$. We first show that
\begin{claim}\label{3}
	$|V(P)\cap C_{1,2}|\leq2$.
\end{claim}
\pf
Suppose to its contrary that $|V(P)\cap C_{1,2}|=3$, then let $V(P)\cap C_{1,2}=\{a,b,c\}$. Since $|V(P)\cap (C_{1,3}\cup C_{2,3})|=1$, we may by symmetry assume that $d\in C_{1,3}$. We prove that
\begin{equation}\label{two vertices}
	\mbox{For any vertex $v_t\in A\setminus\{v_1,v_2\}$, $[v_t,\{a,b,c\}]\ne\emptyset$.}
\end{equation}

Suppose not, and let $v_t\in A\setminus\{v_1,v_2\}$ such that $[v_t,\{a,b,c\}]=\emptyset$. By the definition of $C_{i,j}$ and Claim~\ref{C3}, we have that $\{v_1,v_t\}\cup \{a,b,c,d\}$ is an induced $P_2\cup P_4$, a contradiction. This proves (\ref{two vertices}). 

Recall that $A=\{v_1,\cdots,v_{\omega(G)}\}$ and $\omega(G)\geq5$. By (\ref{two vertices}), each vertex in $\{v_3,v_4,v_5\}$ must have a neighbor in $\{a,b,c\}$. We have that $N_{\{a,b,c\}}(v_i)\cap N_{\{a,b,c\}}(v_j)=\emptyset$ for $i,j\in \{3,4,5\}$ and $i\ne j$, as otherwise, by symmetry, let $a$ be a common neighbor of $v_3$ and $v_4$, but then $\{a,v_3,v_4,v_1\}$ induces a diamond. Therefore, We may suppose that $v_3\sim a$, $v_4\sim b$ and $v_5\sim c$ such that $\{v_3a,v_4b,v_5c\}$ is a matching of $G$. Since $P$ is an induced $P_4$, we have that $P[\{a,b,c\}]$ is either an induced $P_3$ or an $K_2\cup K_1$. By symmetry, if $abc$ is an induced $P_3$, then there is an induced $C_5=v_5v_3abcv_5$. By symmetry, if $a\sim b$ and $[c,\{a,b\}]=\emptyset$, then since $\{a,b,c,d\}$ induces an induced $P_4$, we may suppose $d\sim b$, $d\sim c$ and $d\not\sim a$. There also exists an induced $C_5=v_5v_4bdcv_5$. Both are contradictions. This proves Claim~\ref{3}. \qed

By Claim~\ref{3}, $|V(P)\cap C_{1,2}|\leq2$. If there exists an edge $v_tv_{t'}\in E(G[A\setminus\{v_1,v_2\}])$ such that $[\{v_t,v_{t'}\},V(P)\cap C_{1,2}]=\emptyset$. By the definition of $C_{i,j}$ and Claim~\ref{C3}, $\{v_t,v_{t'}\}\cup V(P)$ induces a $P_2\cup P_4$, a contradiction. Therefore, we have that
\begin{equation}\label{an edge}
	\mbox{For any edge $v_tv_{t'}\in E(G[A\setminus\{v_1,v_2\}])$, $[\{v_t,v_{t'}\},V(P)\cap C_{1,2}]\ne\emptyset$.}
\end{equation}

Since $\omega(G)\geq5$, we have that $|V(P)\cap C_{1,2}|=2$ as otherwise by Claim~\ref{C12}, we can choose such an edge $v_tv_{t'}\in E(G[A\setminus\{v_1,v_2\}])$ such that $[\{v_t,v_{t'}\},V(P)\cap C_{1,2}]=\emptyset$, which contradicts to (\ref{an edge}). Now, we may suppose $V(P)\cap C_{1,2}=\{a,b\}$.

Recall that $A=\{v_1,\cdots,v_{\omega(G)}\}$ and $\omega(G)\geq5$. By (\ref{an edge}), $|N_{A\setminus\{v_1,v_2\}}(a)\cup N_{A\setminus\{v_1,v_2\}}(b)|\geq2$. Moreover, $|N_{A\setminus\{v_1,v_2\}}(a)|\leq1$ and $|N_{A\setminus\{v_1,v_2\}}(b)|\leq1$ according to Claim~\ref{C12}. By symmetry, it implies that $|N_{A\setminus\{v_1,v_2\}}(a)|=|N_{A\setminus\{v_1,v_2\}}(b)|=1$, and $N_{A\setminus\{v_1,v_2\}}(a)\cap N_{A\setminus\{v_1,v_2\}}(b)=\emptyset$. We may by symmetry assume that $a\sim v_3$ and $b\sim v_4$, and now $[a,A\setminus\{v_3\}]=[b,A\setminus\{v_4\}]=\emptyset$.

If $V(P)\cap C_{1,3}=\emptyset$, then $\{v_2,v_5\}\cup V(P)$ is an induced $P_2\cup P_4$ by the definition of $C_{i,j}$ and Claim~\ref{C3}. Hence, $V(P)\cap C_{1,3}\ne\emptyset$ and $V(P)\cap C_{2,3}\ne\emptyset$ by symmetry. We may suppose that $V(P)\cap C_{1,3}=\{c\}$ and $V(P)\cap C_{2,3}=\{d\}$. 

If $a\sim b$, then we may by symmetry assume that $b\sim c$ since $\{a,b,c,d\}$ induces an induced $P_4$. But then, there is an induced $C_5=v_2v_3abcv_2$, a contradiction. Therefore, $a\not\sim b$.

Since $\{a,b,c,d\}$ is an induced $P_4$, $b$ must have a neighbor in $\{c,d\}$. We may by symmetry assume that $b\sim c$. To avoid an induced $C_5=v_1v_4bcdv_1$, $c\not\sim d$. Moreover, since $a\sim c$ leads to an induced $C_5=v_3v_4bcav_3$, we have that $a\not\sim c$. Now, $a\sim d$ and $d\sim b$, i.e., $P=cbda$. But now, there also is an induced $C_5=v_1v_2cbdv_1$, a contradiction. 

It completes the proof of Theorem~\ref{C5}.\qed

\end{document}